\documentclass[11pt,a4paper]{article}

\usepackage{amsmath,amssymb,amsfonts,graphicx,hyperref}
\usepackage[margin=1in]{geometry}

\usepackage{amstext}
\usepackage{stmaryrd}
\usepackage{algpseudocode}

\usepackage[normalem]{ulem}
\usepackage{graphicx}

\usepackage{array}
\usepackage{verbatim}
\usepackage{booktabs}
\usepackage{calc}
\usepackage{textcomp}
\usepackage{multirow}

\usepackage[all]{xy}

\usepackage{cancel}
\usepackage{mathtools}
\usepackage{placeins}
\usepackage{xspace}

\usepackage{tikz,tikzscale}
\usetikzlibrary{calc}
\usetikzlibrary{shadings}

\usepackage[mode=multiuser,status=draft]{fixme} 
\fxsetup{envlayout=color}
\fxsetup{innerlayout=inline}
\fxsetup{targetlayout=colorcb}
\FXProvidesTargetLayout{color}

\fxusetheme{color}
\FXRegisterAuthor{ss}{envss}{SS}
\FXRegisterAuthor{mr}{envmr}{MR}
\FXRegisterAuthor{mw}{envmw}{MW}
\FXRegisterAuthor{mk}{envmk}{MK}

\usepackage{subcaption}
\usepackage[format=plain,indention=.5cm]{caption}
\usepackage{pgfplots}

\newcommand{\Dd}{\text{d}}

\usepackage{bm,upgreek}

\newcommand{\bfm}[1]{{\rm\bf #1}}

\newcommand{\ds}{\,\mathrm{d}s}
\newcommand{\dx}{\,\mathrm{d}x}
\usepackage{graphicx} 
\usepackage[linesnumbered,vlined,commentsnumbered,ruled]{algorithm2e}

\makeatletter
\providecommand\@citeb{}%

\@ifundefined{citation}{%
}{%
    \let\orig@citation\citation
    \def\citation#1{%
        \def\@tempa{#1}%
        \ifx\@tempa\@empty
        \else
        \orig@citation{#1}%
        \fi
    }%
}
\makeatother
\makeatletter
\let\orig@writefile\@writefile
\def\@writefile#1#2{%
    \edef\@tempa{#1}%
    \edef\@tempb{loc}%
    \ifx\@tempa\@tempb
    \else
    \orig@writefile{#1}{#2}%
    \fi
}
\makeatother

\usepackage[final,commandnameprefix=always]{changes} 

\usepackage{xcolor}

\title{Matrix-Free Ghost Penalty Evaluation via Tensor Product Factorization}
\author{Michał Wichrowski\footnote{ Interdisciplinary Center for Scientific Computing,
        Heidelberg University,
        Heidelberg, Germany, mt.wichrowsk@uw.edu.pl} }
\date{}

\begin{document}
\maketitle

\begin{abstract}
    We present a matrix-free approach for implementing ghost penalty stabilization in Cut Finite Element Methods
    (CutFEM). While matrix-free methods for CutFEM have been developed, the efficient evaluation of
    high-order, face-based ghost penalties remains a significant challenge, which this work addresses. By
    exploiting the
    tensor-product structure of the ghost penalty operator, we reduce its evaluation to a
    series of one-dimensional matrix-vector products using precomputed 1D matrices, avoiding the need to evaluate
    high-order derivatives directly. This approach achieves $O(k^{d+1})$ complexity for elements of degree $k$ in $d$
    dimensions, significantly reducing implementation effort while maintaining accuracy.
    The derivation relies on the fact that the cells are aligned with the coordinate axes.
    The method is implemented within the \texttt{deal.II} library.
\end{abstract}

\textbf{Keywords:}  CutFEM, Ghost penalty, Matrix-free, Tensor Product, High-order Finite Elements

\section{Introduction}
Cut finite element methods (CutFEM) have emerged as a versatile approach for solving partial differential equations
(PDEs) on complex geometries \cite{burman2015cutfem, burman2014fictitious, burman2016full, burman2022cutfem,
    cerroni2019numerical}. These methods allow for the discretization of the domain using a mesh that is independent of the
geometry, which simplifies mesh generation and enables the simulation of problems with evolving interfaces or moving
boundaries. However, CutFEM can suffer from ill-conditioning when the interface intersects elements in a way that
creates cut cells with an arbitrarily small volume fraction inside the domain. This can lead to a poorly conditioned
system matrix and adversely affect solver performance. To address this challenge, various stabilization methods have
been developed to ensure mathematical well-posedness regardless of how the geometry intersects with the
interface~\cite{burman2010ghost, badia2022linking}. A popular choice is to employ a ghost penalty approach to penalize
jumps in the solution or its derivatives across element faces, providing stability especially for high-order finite
element spaces~\cite{larson2020stabilization, hansbo2017cut, claus2019cutfem}.

While the direct implementation of face-based ghost penalty can be vexing due to the evaluation of high-order
derivatives, alternatives such as volume-based methods exist that avoid this issue~\cite{burman2010ghost,
    preuss2018higher, bergbauer2024high}. However, the face-based approach is particularly amenable to the tensor-product
factorization presented in this work. By expressing the ghost penalty term as a sum of Kronecker products of
one-dimensional mass and penalty matrices, we can efficiently compute its action on a vector without explicitly
assembling the global matrix. This approach is particularly well-suited for high-order finite element methods, where
the cost of assembling and storing the global matrix can be prohibitively high.

In traditional finite element methods, constructing and storing the stiffness matrix becomes increasingly expensive as
the order of the finite element approximation grows. The fill ratio increases rapidly, leading to significant memory
consumption and data transfer bottlenecks. Matrix-free methods circumvent these limitations by computing the action of
the matrix on a vector on-the-fly, offering substantial performance improvements for high-order finite element
spaces~\cite{kolev2021efficient, kronbichler2012generic, kronbichler2019multigrid, bergbauer2024high}. Recent
works~\cite{davydov2020algorithms, davydov2020matrix, moxey2020efficient, clevenger2020flexible} have demonstrated the
effectiveness of matrix-free methods across various finite element applications~\cite{wichrowski2025finteStrain,
    wichrowski2023exploiting}.

Matrix-free~\cite{kronbichler2012generic} methods have gained popularity in recent years as a way to reduce the
computational cost and memory footprint of finite element simulations. These methods avoid the explicit assembly and
storage of the stiffness matrix, instead computing the action of the matrix on a vector on-the-fly. This can lead to
significant performance improvements~\cite{fischer2020scalability,kronbichler2019fast, kronbichler2019multigrid},
especially for large-scale problems or high-order finite element spaces. The presented method is implemented within the
\texttt{deal.II} library \cite{dealii2019design}, building upon its flexible and efficient framework for matrix-free
methods~\cite{kronbichler2012generic, kronbichler2019multigrid} and the infrastructure developed by Bergbauer
\cite{bergbauer2024high} for cut cells. The implementation of the ghost penalty, however, is novel.

Ghost penalty methods~\cite{burman2010ghost} add stabilization terms to the variational formulation of the PDE,
penalizing jumps in the solution or its derivatives across element faces. Several works have explored the use of ghost
penalty methods in the context of CutFEM. Burman et al. \cite{burman2015cutfem, burman2014fictitious, burman2016full}
introduced various forms of the ghost penalty method and established its mathematical foundations, showing that it
leads to optimal convergence rates independent of how the boundary intersects the mesh. Larson and Zahedi
\cite{larson2020stabilization} analyzed high-order ghost penalty stabilization for CutFEM, proving its stability and
convergence properties for arbitrary polynomial orders. An alternative for the ghost penalty method is using
macro-patch stabilization~\cite{gurkan2019stabilized, preuss2018higher}; that avoids evaluation of derivatives.
However, most ghost penalty approaches, including face-based methods like the one in this manuscript, can be
susceptible to locking depending on the definition of element patches~\cite{badia2022linking, burman2022design}. Trace
Finite Element Method (TraceFEM)~\cite{lehrenfeld2018stabilized, heister2024adaptive, jankuhn2021trace} also relies on
stabilization techniques analogous to the ghost penalty for solving PDEs on surfaces~\cite{heister2024adaptive}.

Badia et al. \cite{badia2022linking} linked strong and weak stabilization methods, discussed ghost penalty locking,
designed locking-free ghost penalties, and compared stabilization schemes. Hansbo et al. \cite{hansbo2017cut} applied
ghost penalty stabilization to cut isogeometric analysis, demonstrating its effectiveness for high-order spline spaces.
Claus and Kerfriden~\cite{claus2019cutfem} used ghost penalty in the context of fluid-structure interaction problems,
where it proved crucial for handling moving interfaces.

More recent developments include the work by Burman et al. \cite{burman2022cutfem} on robust preconditioners for
ghost-penalty-stabilized systems, and {Sticko and Kreiss}~\cite{sticko2016stabilized} on efficient implementation
strategies. {Schoeder}~et al.~\cite{schoeder2020high} presented a high-order accurate cut finite element method with
ghost penalty stabilization for acoustics. Burman et al. \cite{burman2022design} analyzed the design of ghost penalty
operators for various applications, providing guidelines for parameter selection and stability analysis.
In~\cite{cui2025multigrid} a patch smoother for multigrid methods applied to CutFEM yielded promising results.

In~\cite{frachon2024divergence} the authors present a divergence-free cut finite element method for the Stokes
equations, where a special form of the ghost penalty stabilization is used to preserve divergence-free condition while
providing stability. While in this work we focus on the standard ghost penalty stabilization, we expect the presented
method to be extendable to cover this technique too.

In the context of CutFEM, matrix-free evaluation of cell contributions has been explored in~\cite{bergbauer2024high}.
However, applying ghost penalties for high-order methods using their approach would require evaluating higher-order
normal derivatives. To avoid this, volume-based stabilization is often used, though it may lead to locking
issues~\cite{burman2022design}. The method presented in this paper overcomes this issue by demonstrating that
evaluating face-wise ghost penalties only requires precomputed 1D matrices, which can be easily computed regardless of
the method's order.

Together with the benefits of matrix-free evaluation of the finite element operators, comes a challenge of solving the
corresponding linear systems. The matrix-free approach is particularly well-suited for iterative solvers, where the
performance is influenced by the choice of preconditioners. Several preconditioning approaches have been developed for
CutFEM. The work~\cite{gross2021optimal, gross2023analysis} proposed splitting the finite element space into two
subspaces: one spanned by nodal basis functions associated with nodes on the boundary of the fictitious domain, and
another spanned by the remaining nodal basis functions. This decomposition allows for effective preconditioning of the
linear system. A complete matrix-free method for CutFEM was developed by Bergbauer et al.~\cite{bergbauer2024high},
where the authors used projection into an uncut problem to define a preconditioner. Then, a multigrid method was
applied. In this paper, we focus on the matrix-free evaluation of the CutFEM operator with special emphasis on the
ghost penalty operator only, leaving the choice of preconditioners for future work.

Another unfitted approach is the Shifted Boundary Method (SBM)~\cite{main2018shifted, atallah2020second}, which avoids
complex quadrature on cut cells by transferring boundary conditions to a surrogate boundary. This simplifies
implementation while avoiding the need for any stabilization and has been successfully applied to various
problems~\cite{atallah2021shifted}. Recent developments include high-performance matrix-free implementations and
specialized geometric multigrid preconditioners~\cite{wichrowski2025matrix, wichrowski2025geometric,
    wichrowski2026sbmMGCG}. Aggregated Finite Element Methods, merges degrees of freedom in cut cells to ensure stability
without explicit penalty terms~\cite{badia2018mixed, badia2021aggregated}.

This paper presents a matrix-free method for evaluating ghost penalty terms in CutFEM. In
Subsection~\ref{sec:tensor_product}, which forms the core contribution of this work, we demonstrate how the
tensor-product structure of the ghost penalty operator can be exploited for a matrix-free implementation. We show that
the operator is factored into products of one-dimensional operators, which dramatically reduces both the computational
complexity and implementation effort. The key idea is to express the ghost penalty operator in terms of one-dimensional
mass and penalty matrices, which can be precomputed and stored. This allows for an efficient evaluation of the operator
on-the-fly, avoiding the need for explicit assembly of the global matrix.

Key components developed for this work have been integrated into the main \texttt{deal.II} library and are included in
the 9.7 release~\cite{dealII97} This includes the assembly of the 1D ghost penalty stabilization term, the numbering of
degrees of freedom on two adjacent cells, and support for additional ghost cells in the matrix-free
framework\footnote{The source code used for this paper is available at
    \url{https://github.com/mwichro/TensorGhostPenalty}}..

We begin by presenting the CutFEM discretization in Section~\ref{sec:cutFEM}, introducing necessary notation and
concepts. In Section~\ref{sec:numerical_results}, we validate our approach through numerical experiments, demonstrating
optimal convergence rates and computational efficiency. We analyze the performance of the method and its scalability
with respect to polynomial degree and mesh refinement. Finally, Section~\ref{sec:conclusion} summarizes our findings
and discusses potential extensions.

\section{Cut FEM Discretization and Ghost Penalty Stabilization}
\label{sec:cutFEM}
Consider a bounded Lipschitz $d$-dimensional domain $\Omega \subset \mathbb{R}^d$ where we aim to solve the Poisson
problem:
\begin{equation}
    - \Delta u = f \quad \text{in } \Omega,
    \quad u = 0 \quad \text{on } \partial\Omega.
\end{equation}
Traditional FEM partitions $\Omega$ into elements, but CutFEM allows arbitrary intersections with the domain
boundary.
We introduce a Cartesian triangulation \( \mathcal{T}_h \) consisting of square (2D) or cubic (3D) elements of size
$h$. On this mesh, we define a finite element space $\mathbb{V}_h$ using Lagrange polynomial elements of order \( k
\).

To define the discrete domain $\Omega_h$, we introduce a level set function $\phi$ that is positive inside $\Omega$ and
negative outside. The domain $\Omega$ is then discretized using a mesh of cells and faces, where each cell is either
fully inside or at least intersected by the domain.

Our problem is to find the solution $u$ in the space $\mathbb{V}_h$ that satisfies the weak form:
\begin{equation}
    \int_{\Omega} \nabla u \cdot \nabla v \dx = \int_{\Omega} f v \dx \quad \forall v \in \mathbb{V}_h.
\end{equation}
\subsection{Boundary Conditions}
In a classical approach, one requires both the solution and test function to vanish on the Dirichlet boundary. This is
not the case in CutFEM, where it is not possible to enforce this condition directly due to the non-conforming nature of
the mesh. Instead, we weakly enforce the Dirichlet boundary conditions by employing
Nitsche's~\cite{nitsche1971variationsprinzip, warburton2003constants} method. This approach adds a term to the
variational form, defined as
\begin{equation}
    \int_{\partial\Omega} \left( \frac{\partial u}{\partial n} v + \frac{\partial v}{\partial n} u +
    \frac{\gamma_D}{h} u v
    \right) \ds,
\end{equation}
where \( \gamma_D \) is a penalty parameter that provides stability of the method; we use \( \gamma_D = 30 k (k+1) \).

\subsection{Ghost Penalty Stabilization}
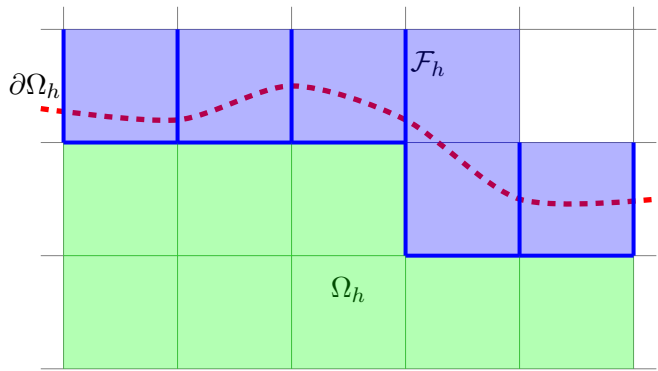
\begin{figure}[!ht]
    \centering
    \begin{tikzpicture}[scale=1.2]
        \draw[step=1cm, gray, very thin] (-0.2,0) grid (5.2,3.2);

        \node at (-0.35,2.5) {$\partial\Omega_h$};
        \draw[line width=2pt, red, dashed] plot [smooth, tension=0.5] coordinates {(-0.2,2.3) (1,2.2) (2,2.5) (3,2.2)
                (4,1.5)
                (5.2,1.5)};

        \node at (2.5,0.7) {$\Omega_h$};
        \node at (3.2,2.7) {$\mathcal{F}_h$};

        \foreach \x/\y in {0/2, 1/2, 2/2, 3/1, 3/2, 4/1}
            {
                \fill[blue, opacity=0.3] (\x,\y) rectangle (\x+1,\y+1);
            }
        \foreach \x/\y in {0/0, 0/1,2/1, 1/0, 1/1, 2/0, 3/0, 4/0}
            {
                \fill[green, opacity=0.3] (\x,\y) rectangle (\x+1,\y+1);
            }

        \foreach \x/\y in {0/2, 1/2, 2/2, 2/1, 3/1 , -1/2, 4/1	} {
                \draw[blue, line width=1.5pt] (\x+1,\y) -- (\x+1,\y+1);
            };

        \foreach \x/\y in {0/1, 1/1, 2/1		/	,	3/0, 4/0, 3/1 } {
                \draw[blue, line width=1.5pt] (\x,\y+1) -- (\x+1,\y+1);
            };

    \end{tikzpicture}
    \caption{Illustration of a computational domain discretization for CutFEM: The curved domain boundary
        $\partial\Omega_h$ (red dashed line) intersects a Cartesian mesh, creating cut cells (blue) and interior cells
        (green).
        Ghost penalty stabilization is applied across the faces $\mathcal{F}_h$ (thick blue lines) between cut cells or
        between
        cut and interior cells to ensure numerical stability regardless of the boundary position.}
    \label{fig:intersected_cells}
\end{figure}

The ghost penalty term
\begin{equation}
    g_h(u_h,v_h) = \gamma_A \sum_{F \in \mathcal{F}_h} g_F(u_h, v_h)
\end{equation}
penalizes discontinuities by adding stabilization contributions over faces of elements cut by the
boundary, as depicted in Figure~\ref{fig:intersected_cells}. We denote the set of all such faces as $\mathcal{F}_h$.
This operator penalizes jumps of derivatives over faces, and handwavingly speaking, it glues pieces of polynomial
functions
on two adjacent cells together into a single polynomial function. The penalty term is defined as:
\begin{equation}
    g_F(u_h, v_h) = \gamma_A \sum_{k=0}^p \left(\frac{h_F^{2k-1}}{k!^2}[\partial_n^k u_h], [\partial_n^k v_h]
    \right)_F
\end{equation}
with $h_F$ being the diameter of the face, $\partial_n^k$ the normal derivative of order $k$, and
$[\cdot]_F$ the jump
across the face. The parameter $\gamma_A$ is chosen to balance the penalty and the original bilinear form. This
provides stability of the system, independent of the location of the
interface~\cite{burman2014fictitious}.

There exists an alternative approach to the ghost penalty stabilization, namely macro patch stabilization, where the
system is stabilized by integrating the difference between the function and the extension of the function from the
neighboring element. The derivation follows the same idea as presented here, and it exhibits the identical tensor
product structure.

The issue with the ghost penalty is that it requires evaluation of high-order derivatives across faces, which can be a
challenging task, especially in a matrix-free setting. This involves handling higher-order derivatives of the mapping,
which in the case of arbitrary meshes can be quite complex. In CutFEM, however, we can simplify this by assuming that
the mesh is Cartesian, where the shape of each cell is defined by a size in each dimension. Moreover, the local matrix
of the ghost penalty term is identical for each face; thus, it can be precomputed and stored in memory. The evaluation
cost of the ghost penalty term can be reduced by exploiting the tensor product structure of the shape functions.

\subsection{Weak Formulation}
The weak form of the problem is to find $u\in \mathbb{V}_h$ such that
\begin{equation}
    \int_{\Omega} \nabla u \cdot \nabla v \dx + \int_{\partial\Omega} \left( \frac{\partial u}{\partial n} v +
    \frac{\partial v}{\partial n} u + \frac{\gamma_D}{h} u v \right) \ds + \sum_{F \in
        \mathcal{F}_h} g_F(u,v)= \int_{\Omega} f v
    \dx \quad\quad \forall v \in
    \mathbb{V}_h.
\end{equation}

\subsection{Exploiting the Tensor Product Structure}
\label{sec:tensor_product}
For simplicity of presentation, we consider a 3D case with a hexahedral element; however, the reasoning can be easily
transferred to the 2D case with quadrilateral elements. The derivation relies on the fact that the cells are
aligned with the coordinate axes. Consider a family of basis functions that are products of
one-dimensional functions. Let the shape functions in
each direction be denoted by \( \phi_i(x_1) \), where \( x_1 \) is the
reference coordinate in one dimension.
We introduce the tensor product numbering, where each degree of freedom in a higher-dimensional element (e.g., a
quadrilateral or hexahedral element) is associated with a combination of 1D indices, one for each spatial direction.
These indices are combined into a multi-index $(i_1, i_2,     i_3)$,
allowing us to efficiently reference and compute values for multi-dimensional functions using one-dimensional basis
functions
\begin{equation}
    \phi_{i_1, i_2, i_3}( x_1, x_2, x_3 ) = \phi_{i_1}(x_{1}) \cdot \phi_{i_2}(x_{2}) \cdot \phi_{i_3}(x_{3}).
    \label{eq:tensor_product}
\end{equation}

The values of function \( u \) at a given point are expressed as a linear combination of these basis functions with
coefficients \( \bm{[}u_{i_1 i_2 i_3}\bm{]} \). Here, the multi-index \( (i_1, i_2, i_3 )\) corresponds to the tensor
product indexing of degrees of freedom associated with the element. The value of the function at point \( x \) with
coordinates \( (x_{1}, x_2, x_3) \) is given by
\begin{equation}
    u (x) =u(x_{1}, x_2, x_3) = \sum_{i_1, i_2, i_3} u_{i_1 i_2 i_3} \phi_{i_1}(x_{1}) \cdot \phi_{i_2}(x_{2}) \cdot
    \phi_{i_3}(x_3)
\end{equation}
In practice, coefficients  $u_{i_1 i_2 i_3}$ are stored in a one-dimensional array using lexicographical ordering.
The conversion between the one-dimensional index $i$ and the corresponding multi-index $(i_1, i_2, i_3)$ is
done using a simple formula
\begin{equation}
    i = i_1 + i_2 (N_1) + i_3 (N_1 N_2)
\end{equation}
where \( N_1, N_2, N_3 \) are the number of basis functions in each direction.

\begin{figure}[!ht]
    \centering
    \begin{tikzpicture}[scale=1.7]
        \draw[line width=1pt, dashed] (-2,0) rectangle (0,2) node[midway, below] {$K_1$};
        \foreach \x in {0,1,2,3, 4}
        \foreach \y in {0,1,2 }
        {
        \node[circle,inner sep=2pt,fill=red!70,label={[label distance=2pt]below left:{\small (\x,\y)}}]
        (K1-\x-\y) at (-2+\x*2/2,\y*2/2) {};
        \pgfmathtruncatemacro{\index}{\x + \y * 5}
        \node[above right=1pt,font=\small] at (K1-\x-\y) {\index};
        }

        \draw[line width=1pt, dashed] (0,0) rectangle (2,2) node[midway, below] { $K_2$};

        \draw[line width=2pt, green!70] (0,0) -- (0,2) node[midway, right, below, black] {  $\; \;\; F$};

    \end{tikzpicture}
    \caption{Tensor product numbering of degrees of freedom for two adjacent quadratic elements $K_1$ and $K_2$  in 2D,
        sharing a face $F$. The numbers indicate the lexicographical ordering of the DoFs within each cell as well
        as the corresponding multi-indices. The figure shows a 2D case for simplicity of presentation.}
    \label{fig:tensor_product_numbering}
\end{figure}
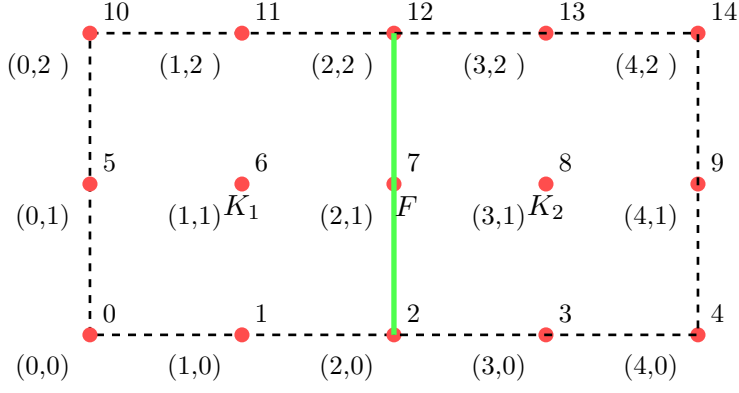

We consider a pair of cells \( K_1 \) and \( K_2 \) that share a face \( F \) orthogonal to the \( x_1 \) axis, as
illustrated in Figure~\ref{fig:tensor_product_numbering}. Since we consider a Cartesian mesh, the face is a
hyperrectangle. It can be represented as a Cartesian product of intervals $F_i$. Specifically, $F = F_1 \times F_2
    \times F_3$ where $F_1$ contains a single point.

Next, we define the local basis $ \{ \phi_i \}$ by numbering the shape functions lexicographically and observe that
they form a tensor product structure with $ N_1=2k+1$ and $ N_i=k+1$ for $i\neq 1$. We will further elaborate on the
implications of this structure in the context of ghost penalty terms. For simplification we will consider only the term
with \( k \)-th derivative in the ghost penalty term on the face \( F \):
\begin{equation}
    g_{F,k}(u_h, v_h) =
    \left( [\partial_n^k u], [\partial_n^k v]  \right)_F .
\end{equation}
The local penalty matrix $\mathcal{G}_{F,k}$ for the said face is given by
\begin{equation}
    [\mathcal{G}_{F,k} ]_{i,j} =
    \left( [\partial_n^k \phi_j], [\partial_n^k \phi_i] \right)_F.
\end{equation}
We break down indices \( i \) and \( j \) into multi-indices \( (i_1, i_2, i_3) \) and \( (j_1, j_2, j_3)\). By
expanding the shape functions using tensor products of one-dimensional shape functions, as in
Equation~(\ref{eq:tensor_product}) we obtain
\begin{align}
    [\mathcal{G}_{F,k} ]_{i,j} &
    = \left( [\partial_n^k \phi_j], [\partial_n^k \phi_i] \right)_F                             \\
                               & = \Bigl( [\partial_n^k \phi_{i_1}] \cdot \phi_{i_2}\phi_{i_3},
    [\partial_n^k \phi_{j_1}] \cdot \phi_{j_2}\phi_{j_3} \Bigr)_F                               \\
                               & =  \left([\partial_n^k \phi_{i_1}] \cdot [\partial_n^k
        \phi_{j_1}]\right)_{F_1} \cdot
    \left( \phi_{i_2}\phi_{i_3} , \phi_{j_2}\phi_{j_3} \right)_{F_2\times F_3}
\end{align}
The part $  [\partial_n^k \phi_{i_1}]_{F_1} \cdot [\partial_n^k     \phi_{j_1}]_{F_1} $
is   1D ghost penalty term corresponding to the  jump of the normal derivative of the 1D shape functions
\begin{equation}
    [G^h_k]_{i_1,j_1} = \left[\partial_n^k \phi_{i_1}\right]_{F_1} \, \left[\partial_n^k
        \phi_{j_1}\right]_{F_1}.
\end{equation}
The remaining part is the scalar product of the shape functions over the face can be broken down to 1D integrals
\begin{equation}
    \Bigl( \phi_{i_2}\phi_{i_3} , \phi_{j_2}\phi_{j_3} \Bigr)_{F_2\times F_3}
    = \Bigl( \phi_{i_2} , \phi_{j_2} \Bigr)_{F_2} \cdot \Bigl( \phi_{i_3} , \phi_{j_3} \Bigr)_{F_3}.
\end{equation}
where each corresponds to the mass matrix in the other 1D
\begin{equation}
    M^{2,h }_{i_2,j_2} = \Bigl( \phi_{i_2} , \phi_{j_2} \Bigr)_{F_2}, \quad
    M^{3,h}_{i_3,j_3} = \Bigl( \phi_{i_3} , \phi_{j_3} \Bigr)_{F_3}.
\end{equation}
It is easy to see that all of those mass matrices are identical. Note that the number of basis functions in directions
parallel to the face is $k+1$ and the number of basis functions in the direction orthogonal to the face is
2k+1.

A simple calculation shows that in this setting the ghost penalty matrix $\mathcal{G}_{F,k}$ can be expressed as the
following Kronecker product
\begin{equation}
    \mathcal{G}_{F,k} = M^{h} \otimes M^{h} \otimes G^h_k .
\end{equation}
The matrices are precomputed using shape functions $\hat{\phi}_\circ $ on the reference  1D cell $\hat{K}$ with $h=1$
and then adjusted to the cell geometry. The face of the 1D cell is a point denoted by $\hat{F}$, and the 1D ghost
penalty matrix
is given by
\begin{align}
    [G^h_k]_{i_1,j_1} & = \left[\partial_n^k \phi_{i_1}\right]_{F}\cdot \left[\partial_n^k \phi_{j_1}\right]_{F}
    =  \left[{1\over h^k} \partial_n^k \phi^{\text{ref}}_{i_1}\right] _{\hat{F}} \cdot \left[{1\over h^k}
        \partial_n^k \phi^\text{ref}_{j_1}\right]_{\hat{F}}						     \nonumber
    \\
                      & = {1\over h^{2k}} \left[\partial_n^k \phi^{\text{ref}}_{i_1}\right]_{\hat{F}}\cdot
    \left[\partial_n^k \phi^{\text{ref}}_{j_1}\right]_{\hat{F}}= {1\over h^{2k}} [G^1_k]_{i_1,j_1}.
\end{align}
The mass matrix is scaled by the cell diameter $h$ to obtain the mass matrix $M^h$ for the cell $K$:
\begin{equation}
    [M^h]_{i_1,j_1} = \left[\phi_{i_1}, \phi_{j_1}\right]_K =h \left(\phi^{\text{ref}}_{i_1},
    \phi^{\text{ref}}_{j_1}\right)_{\hat{K}} = h[M^1]_{i_1,j_1}.
\end{equation}
Finally, the complete ghost penalty matrix is then given by the following sum
\begin{align}
    \mathcal{G}_F & = \gamma_A \sum_{k=0}^p \frac{h_F^{2k-1}}{k!^2} \left( M^h \otimes M^h \otimes G^h_k
    \right)  \nonumber
    \\
                  & = \gamma_A \;	M^h \otimes M^h \otimes \sum_{k=0}^p \left(
    \frac{h^{2k-1}}{k!^2}
    G^h_k \right) \nonumber
    \\
                  & = \gamma_A \; M^h \otimes M^h \otimes	\sum_{k=0}^p \left( {  h^{-1} k!^{-2}}
    G^1_k
    \right)  \nonumber
    \\
                  & = \gamma_A \; h M^1 \otimes h M^1 \otimes	\sum_{k=0}^p \left( { h^{-1} k!^{-2}}
    G^1_k
    \right).
    \label{eq:ghostPenaltyMatrixRev}
\end{align}
We store the precomputed mass matrix $h M^1$ and the penalty matrix
\begin{equation}
    G^1 = \sum_{k=0}^p \left( { h^{-1} k!^{-2}} G^1_k \right).
\end{equation}

A similar approach can be used to compute the ghost penalty matrix for the faces orthogonal to the other axes. The
resulting formula for the ghost penalty matrix on face orthogonal to the $x_2$ and $x_3$ axis ($F^2$ and $F^3$
respectively) are given by:
\[
    \mathcal{G}_{F^2} = h M^1  \otimes G^h \otimes h M^1, \quad \mathcal{G}_{F^3} = G^h  \otimes  h M^1 \otimes h M^1 .
\]


\subsection{Matrix-Free Operator Application}

Applying the stabilized CutFEM operator to a vector $ w= \mathcal{A}u$ can be done by looping over all cells and faces
and faces with ghost penalty terms and accumulating their contributions. We will first discuss the application of the
ghost penalty part of the operator as it is the main focus of this work.

\subsubsection{Application of the Ghost Penalty Operator}

The ghost penalty operator is applied to a vector $u$ by computing the contributions of the ghost penalty term on each
face. The contribution of the face $F$ to the vector $w$ is given by
\begin{equation}
    w_F =  \gamma_A  h^{2} M^h \otimes M^h \otimes G^1	\cdot u.
\end{equation}
The application of the ghost penalty operator is decomposed into a sequence of matrix-vector products, which are
computed in three steps:
\begin{align}
    w_F^1 & =  \gamma_A \;  I \otimes I \otimes G^1 \cdot u  ,       \\
    w_F^2 & =  \gamma_A \;  I \otimes h M^1 \otimes I \cdot w_F^{1}, \\
    w_F   & =  \gamma_A \; h M^1 \otimes I \otimes I \cdot w_F^{2}.
\end{align}
The first step computes the contribution of the ghost penalty term on the face $F$ to the vector $u$. The
second
and third steps apply the mass matrix $M^h$ in the directions parallel to the face.
A simple calculation shows that the cost per face of applying the 1D ghost penalty
operator is $(k+1)^{d-1}(2k+1)^2$ operations, and applying each of the two mass matrices costs $(k+1)^d(2k+1)$
operations, leading to a total cost of $(k+1)^{d-1}(2k+1)^2 + 2(k+1)^d(2k+1)$ operations, which is $O(k^{d+1})$.

\subsubsection{Treatment of Internal Cells}
The treatment of interior cells follows the standard matrix-free approach and is included here for completeness. For
each cell $K \in \mathcal{T}$, the $i$-th element of the local contribution ${w}_{Ki}$ is computed by
\begin{equation}
    w_{Ki}=\int\limits_K  \chreplaced{\nabla_x}{\nabla} u_k   \cdot \chreplaced{\nabla_x}{\nabla} {\phi }_i \;	\Dd x
    = \sum_{q(K)}	\chreplaced{\nabla_x}{\nabla} u_k   \cdot  \chreplaced{\nabla_x}{\nabla}\phi_i \, J(x_q) \,
    \omega_q.
\end{equation}
where $\phi_i$ is the $i$-th function of the local basis. The sum is taken over the $q(K)$ quadrature points $x_q$,
with $J(x_q)$ being the Jacobian determinant
at  $q$th quadrature point, and $\omega_q$ the quadrature weight. \chadded{Here and throughout the paper, we use $\nabla_x$ to denote the gradient with respect to the physical coordinates $x$, whereas $\nabla_\xi$ denotes the gradient with respect to the coordinates $\xi$ of the reference cell.}

The evaluation of this term is split into three stages. First, the \textit{evaluation} step: for each quadrature point
$x_q$, \chadded{we compute the gradient of the solution in reference coordinates $\chreplaced{\nabla_\xi}{\nabla}
        u_K(x_q)$, transform it to the physical coordinates using the inverse-transpose of the Jacobian matrix $\bfm{J}^{-T}$
    to obtain ${\nabla_x} u_K(x_q)$,} and store it. Next, in the \textit{quadrature loop} step, we prepare the data for
integration by queuing the computed gradients $\chreplaced{\nabla_x}{\nabla} u_K(x_q)$. Finally, in the
\textit{integration} step, we compute the local contributions $w_{Ki}$ by performing the integration, which involves
contracting the stored gradients with test function gradients $\chreplaced{\nabla_x}{\nabla}\phi_i(x_q)=
    \chadded{\bfm{J}^{-T} \nabla_\xi \phi_i(x_q) }$, Jacobian determinants $J(x_q)$, and quadrature weights $\omega_q$.
Both stages are performed efficiently using sum factorization, where the integration step corresponds to the transpose
of the evaluation operation~\cite{kopriva2009spectral, deville2002high}.

Matrix-free methods typically use sum factorization~\cite{kronbichler2019fast,fischer2020scalability} to bring down the
cost of evaluating gradients \( \chreplaced{\nabla_x}{\nabla} u_K( {x}_q) \) at quadrature points. For integration, we
use a tensor-product quadrature formula where points are indexed with a multi-index \( (q_1, q_2, q_3 )\) and each
quadrature point \( x_{q} = (x_{q_1}, x_{q_2}, x_{q_3} ) \) is represented by a vector of corresponding points of the
one-dimensional quadrature formula. For a finite element of degree $k$, we use $(k+1)$ quadrature points in each
direction.

We define the matrix \(\bm{[}\phi_i(x_{q_j})\bm{]} \), which contains the values of the 1D shape functions \( \phi_i \)
at the 1D quadrature point \( x_{q_j} \). A straightforward calculation shows that the set of function values at
quadrature points $u_{K}^q $ is a result of the following operation:
\begin{equation}
    u_K^q = \left( \bm{[}\phi_i(x_{q_1})\bm{]} \otimes \bm{[}\phi_i(x_{q_2})\bm{]} \otimes \bm{[}\phi_i(x_{q_3})
    \bm{]} \right)u_K.
    \label{eq:u_q}
\end{equation}

Sum factorization exploits the separable structure of tensor-product elements to reduce the computational cost of
evaluation of the formula~(\ref{eq:u_q}). Instead of evaluating the multidimensional, the evaluation is split into a
series of one-dimensional operations:
\begin{align}
    u_{K}^1 & = \left( \bm{I} \otimes \bm{I} \otimes \bm{[}\phi_i(x_{q_1})\bm{]} \right) u_K  ,  \\
    u_{K}^2 & = \left( \bm{I} \otimes \bm{[}\phi_i(x_{q_2})\bm{]} \otimes \bm{I} \right) u_K^1 , \\
    u_{K}^q & = \left( \bm{[}\phi_i(x_{q_3})\bm{]} \otimes \bm{I} \otimes \bm{I} \right) u_K^2 .
\end{align}
Each of the operations above is equivalent to applying the matrix $\bm{[}\phi_i(x_q)\bm{]}$ to the
corresponding subranges of the array $\bm{[}u_{i_1,i_2,i_3}\bm{]}$. Thus, we avoid the need for a full
multidimensional evaluation, thereby
reducing the complexity to $\mathcal{O}((k+1)^{d+1}) = \mathcal{O}((k+1)     N_c)$, where $N_c = (k+1)^d$ is the number
of degrees of freedom per element. The benefits are especially visible for higher-order elements, transforming what
would otherwise be quadratic growth in computational cost into an almost linear complexity with respect to the degrees
of freedom. Note that the cost of the cell evaluation grows at the same rate as the cost of evaluating the ghost
penalty.

\subsection{Treatment of Intersected Cells}

For cells intersected by the boundary, the integration is performed only over the part of the cell that lies inside the
domain $\Omega$. The quadrature formula for such cells is generated according to a method described in
\cite{saye2015high}, and implemented in \texttt{deal.II}~\cite{dealii2019design}. The generated quadrature rule does
not follow the tensor product structure of the quadrature rule and as the results sum factorization cannot be used. The
quadrature points are generated using the mapping of the cell to the reference cell, with the quadrature points and
weights adapted to the actual geometry of the intersection. For each intersected cell, we store:
\begin{itemize}
    \item locations of quadrature points $x_q$ within the cell, \item{ inverse of the Jacobian matrix $\bfm{J}^{-T}$ of the
          reference-to-real cell mapping,}
    \item integration weights $\omega_q$ accounting for the cut geometry,
    \item normal vectors at quadrature points on the boundary.
\end{itemize}

The computation of the local contribution follows the same general pattern as for interior cells, but with
modifications to account for the cut geometry. The weak form of the equation is approximated as
\begin{equation}
    w_{Ki}=\int\limits_{K \cap \Omega}	\chreplaced{\nabla_x}{\nabla} u_K   \cdot \chreplaced{\nabla_x}{\nabla} {\phi }_i \;	\Dd x
    \approx \sum_{q(K \cap \Omega)}	\chreplaced{\nabla_x}{\nabla} u_K   \cdot  \chreplaced{\nabla_x}{\nabla}\phi_i \, J(x_q) \, \omega_q.
\end{equation}

The evaluation of the local contribution proceeds in several steps, ordered to optimize memory access by reading the
degrees of freedom for a cell only once.
\begin{enumerate}
    \item \textbf{Evaluation of Solution Gradients:} Evaluate the solution values at each quadrature point $x_I$
          inside the
          domain:
          \begin{align}
              \chreplaced{\nabla_x}{\nabla} u_K(x_I) & = \sum_i u_{K,i} \bfm{J}(x_I)^{-T}	\chreplaced{\nabla_\xi}{\nabla} \phi_i(x_I).
          \end{align}
          Note: the mapping of gradients from the reference cell to the physical cell uses the
          inverse-transpose
          of the Jacobian matrix. In formulas below we therefore apply
          $\nabla_x = \bfm{J}^{-T}\,\nabla_{\xi}$, where $\bfm{J}$ denotes the Jacobian matrix of the
          reference-to-physical mapping and $J$ its determinant (used in quadrature weights).
    \item \textbf{Evaluation of Gradients and Values at the Boundary:} Evaluate the solution values at each
          quadrature point $ x_S$ at the boundary:
          \begin{align}
              u_K(x_S)                               & = \sum_i u_{K,i}  \phi_i(x_S)                                                   \\
              \chreplaced{\nabla_x}{\nabla} u_K(x_S) & = \sum_i u_{K,i} \bfm{J}^{-T}(x_S)	\chreplaced{\nabla_\xi}{\nabla} \phi_i(x_S).
          \end{align}
    \item \textbf{Integration inside the cell:} Integrate using the precomputed weights specific to the part of the
          cell inside the domain:
          \begin{equation}
              w^{\Omega}_{Ki} = \sum_{q(K \cap \Omega)} \chreplaced{\nabla_x}{\nabla} u_K(x_q) \cdot \chreplaced{\nabla_x}{\nabla}\phi_i(x_q) \cdot J(x_q) \cdot
              \omega_q.
          \end{equation}

    \item \textbf{Integration on the boundary:} Integrate using the precomputed weights specific to the surface:
          \begin{equation}
              w^{\partial\Omega}_{Ki} \stackrel{+}{\leftarrow} \sum_{q(K \cap \partial\Omega)} \left( \frac{\partial
                  u}{\partial
                  n}(x_q)
              \phi_i(x_q) + \frac{\partial \phi_i}{\partial n}(x_q) u(x_q) + \frac{\gamma_D}{h} u(x_q)
              \phi_i(x_q) \right)
              \cdot
              J(x_q)
              \cdot
              \omega_q.
          \end{equation}
    \item  \textbf{Sum up contributions:}
          \begin{equation}
              w_{Ki}= w^{\Omega}_{Ki} + w^{\partial\Omega}_{Ki}.
          \end{equation}

\end{enumerate}

\subsection{Implementation Details}

As discussed, the operator is applied by looping over all cells and faces, accumulating the contributions of the ghost
penalty and the interior cells. The procedure of applying the operator $\mathcal{A}$ is illustrated by
Algorithm~\ref{alg:mf_generic}. This algorithm outlines the steps for computing $w = \mathcal{A}u$ in a matrix-free
manner, accounting for both interior and intersected cells, as well as the ghost penalty contributions.

As the matrix-vector product is a key performance component of any iterative solver, we aim to utilize the hardware as
efficiently as possible. Hence using vectorized operation is crucial for the performance of the method. On every cell
inside the domain the operation is identical hence vectorization over the cells comes up naturally. The other group of
cells are the intersected cells, where the number quadrature points may vary from cell to cell.

We assign each cell one of two categories: interior cells and intersected cells and next group cell of each category
into batches. While applying the operator we loop over the batches of cells. The interior cells are processed in a
vectorized manner. The intersected cells in each batch are processed one by one and SIMD instructions are used to
vectorize the operations within one cell. Our implementation relies on \texttt{deal.II}~\cite{dealii2019design}, in
particular infrastructure developed in the work by Bergbauer~\cite{bergbauer2024high}.

The degrees of freedom are numbered lexicographically within each cell, as illustrated in
Figure~\ref{fig:tensor_product_numbering}. However, for the ghost penalty operator, we require a lexicographical
ordering across two neighboring cells. To achieve this, we copy the values of local degrees of freedom from both cells
into a single array with the desired ordering. Then the operator is applied to the combined array, and the results are
scattered back to the local arrays.

\chadded{To ensure efficient utilization of parallel hardware, we use distributed computing via MPI through
    infrastructure implemented in \texttt{deal.II}. The cells are partitioned and distributed among the processors, and
    the operator is applied in parallel. For load balancing, we use a similar strategy to~\cite{bergbauer2024high} where
    to each cell a weight is assigned. For cells outside the domain the weight is zero, for cells strictly inside the
    domain the weight is one, and for intersected cells the weight is $k^{d-1}$. This estimate comes from the comparison
    of complexity of evaluation in interior cells and intersected cells. The cells are then redistributed among the
    processors (MPI ranks) in such a way that the sum of the weights is balanced across the partition.}

The ghost penalty operator is applied in a separate loop over the relevant faces, which are grouped into vectorization
batches based on their normal direction. In the parallel setting, faces shared between MPI ranks are processed by the
rank with the lowest ID to avoid redundant computations. This simple rule can lead to a slight load imbalance, but its
impact is limited since the ghost penalty evaluation constitutes a relatively small fraction of the total computational
cost, as shown in Figure~\ref{fig:vmult_time_breakdown}.

\chadded{The communication between the processors is minimized
    by using ghosted vectors, which store the values of the solution on the neighboring cells. The values of the solution
    on the ghost cells are updated after each application of the operator. The main difference between the application of
    the operator considered in this work and the standard matrix-free approach is that the ghost penalty operator is
    applied to the faces. This requires extending the data structures inside the matrix-free infrastructure to store
    information about the ghost cells.}

Further optimizations, such as fusing the face loop with the cell loop down to the level of memory transfers, were
considered but not implemented due to their complexity and the diminishing returns as the mesh is refined and the
fraction of stabilized faces decreases.

\begin{algorithm}[!ht]
    \SetKwInOut{Input}{Given}
    \SetKwInOut{Output}{Return}
    \Input{  ${u}$ - current FE solution \\}
    \Output{ ${w}=	\mathcal{A}u$ }
    $w \leftarrow 0$  \tcp*{Initialize destination vector}
    Update ghost values of $u$ \tcp*{Parallel communication}

    \ForEach{element $K \in \Omega^h$, $K \cap \Omega \neq \emptyset$  }{
        \If{$K\subset \Omega $	\tcp*{Handle as interior cell}	}{

            Gather element-local vector values\\ Evaluate gradients at each quadrature point \chadded{$x_q$}: \\ \Indp\chadded{
                $\nabla u_q$}\tcp*{Sum factorization} \Indm \ForEach{quadrature point $q$ on $K$ }{ Queue
                $\chreplaced{\nabla_x}{\nabla} u_q$ for integration\; } Integrate queued gradients: \\ \Indp $w_{Ki} \leftarrow \sum_q
                \chreplaced{\nabla_\xi}{\nabla} {\phi}_i(x_q) \cdot \bfm{J}^{-T}(x_q) \chreplaced{\nabla_x}{\nabla} {u}_K(x_q) \, \,
                \omega_q$ \tcp*{Sum factorization} \Indm Scatter results to ${w}$\\

        }

        \ElseIf{$K \cap \partial \Omega^h \neq 0$  \tcp*{Handle as intersected cell} }{
            Gather element-local vector

            Evaluate gradients and values at surface quadrature point: \\ \Indp $\chreplaced{\nabla_x}{\nabla} {u}_{K}(x_S) \qquad
                {u}_{K} (x_S) $\\ \Indm

            \ForEach{quadrature point $q$ inside  $\Omega$	}{
                Queue $u(x_S)$	and   $ \frac{\partial u}{\partial n}(x_S)	+ \frac{\gamma_D}{h}
                    u(x_S)  $ for integration \\
            }
            Integrate queued gradients:\\
            \Indp
            $w_{Ki} \leftarrow \sum_ {x_S}\chreplaced{\nabla_\xi}{\nabla} {\phi}_i(x_S) \cdot  \bfm{J}^{-T}(x_S) \chreplaced{\nabla_x}{\nabla}
                {u}_K(x_S) \,
                \, \omega(x_S) $
            \\
            \Indm

            Evaluate gradients at each quadrature point inside the domain: \\ \Indp $\chreplaced{\nabla_x}{\nabla} {u}_{K}(x_I) $\\
            \Indm \ForEach{quadrature point $q$ on $\partial \Omega$ }{ Queue $\chreplaced{\nabla_x}{\nabla} u_q$ for integration
                \\ } Integrate queued gradients: \\ \Indp $w_{Ki} \stackrel{+}{\leftarrow} \sum_q \chreplaced{\nabla_\xi}{\nabla}
                {\phi}_i(x_q) \cdot \bfm{J}^{-T}(x_q) \chreplaced{\nabla_x}{\nabla} {u}_K(x_q) \, \, \omega_q$ \\ \Indm Scatter results
            to ${w}$\\ } } \ForEach{face $F \in \mathcal{F}_h^h$ \tcp*{Apply ghost penalty} }{ Gather face-local vector values on
        this face $u_F$\\ Evaluate local ghost penalty operator: \\ \Indp \(w_F \leftarrow \mathcal{G}_F u_F \) \\ \Indm
        Scatter results to ${w}$\\ } Import contributions $w$ \tcp*{Parallel communication} \caption{Matrix-free application of
        the CutFEM operator} \label{alg:mf_generic}
\end{algorithm}

\section{Results}
\label{sec:numerical_results}
To validate the effectiveness of our matrix-free CutFEM implementation with ghost penalty stabilization, we present a
series of numerical experiments. These results demonstrate the method's convergence properties, computational
efficiency, and memory transfer characteristics. As the main difference between the standard matrix-free approach
and the one considered are operations involving the intersected cells (either directly or via ghost penalty operator)
we also demonstrate the impact of the fraction of intersected cells on the overall performance.
The computations were performed on a machine equipped with dual-socket AMD EPYC 7282 processor at 2800 GHz with 16
cores
per socket, the performance measurements were taken using 32 MPI ranks.
The code was compiled using GCC 11.4.0 with the optimization level set to \texttt{-O3}.

\subsection{Computational Efficiency}

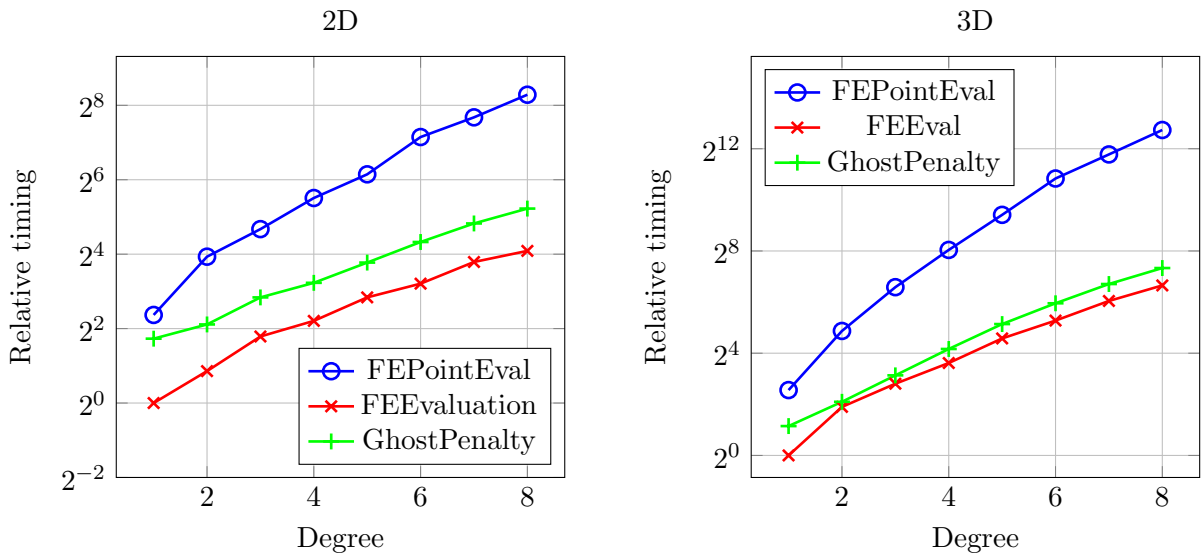
\begin{figure}[!ht]
    \centering
    {\label{fig:single_cell_performance_2d}
        \begin{tikzpicture}
            \begin{axis}[
                    title=2D,
                    xlabel=Degree,
                    ylabel=Relative timing,
                    legend pos=south east ,
                    width=0.47\textwidth,
                    height=0.45\textwidth,
                    grid=both,
                    ymode=log,
                    log basis y=2,
                    ymin =0.25
                ]
                \addplot[color=blue, mark=o, mark size=3pt, line width=1pt] table [x index=0, y
                        index=1, y expr=\thisrowno{1}*4/(7.7400e-02), header=has
                        colnames]
                    {results/one_cell_2d.txt};
                \addlegendentry{\texttt{FEPointEval}}

                \addplot[color=red, mark=x, mark size=3pt, line width=1pt] table [x index=0, y
                        index=2,  y expr=\thisrowno{2}/(7.7400e-02), header=has
                        colnames]
                    {results/one_cell_2d.txt};
                \addlegendentry{\texttt{FEEvaluation}}

                \addplot[color=green, mark=+, mark size=3pt, line width=1pt] table [x index=0, y
                        expr=\thisrowno{3}/(7.7400e-02),
                        header=has
                        colnames]
                    {results/one_cell_2d.txt};
                \addlegendentry{\texttt{GhostPenalty}}

            \end{axis}
        \end{tikzpicture}
    }
    \hfill
    {\label{fig:single_cell_performance_3d}
        \begin{tikzpicture}
            \begin{axis}[
                    title=3D,
                    xlabel=Degree,
                    ylabel=Relative timing,
                    legend pos=north west,
                    width=0.47\textwidth,
                    height=0.45\textwidth,
                    grid=both,
                    ymode=log,
                    log basis y=2,
                    ymin =0.55,
                    ymax = 0.5e5
                ]
                \addplot[color=blue, mark=o, mark size=3pt, line width=1pt] table [x index=0, y
                        index=1, y expr=\thisrowno{1}*4/(2.2180e-01), header=has
                        colnames]
                    {results/one_cell_3d.txt};
                \addlegendentry{\texttt{FEPointEval}}

                \addplot[color=red, mark=x, mark size=3pt, line width=1pt] table [x index=0, y
                        index=2,  y expr=\thisrowno{2}/(2.2180e-01), header=has
                        colnames]
                    {results/one_cell_3d.txt};
                \addlegendentry{\texttt{FEEvaluation}}

                \addplot[color=green, mark=+, mark size=3pt, line width=1pt] table [x index=0, y
                        expr=\thisrowno{3}/(2.2180e-01),
                        header=has
                        colnames]
                    {results/one_cell_3d.txt};
                \addlegendentry{\texttt{GhostPenalty}}

            \end{axis}
        \end{tikzpicture}
    }
    \caption{Relative time per application for different evaluation methods on a single cell. The time is
        normalized by the time of \texttt{FEEvaluation} for $k=1$. In 2D that is 0.077 $\mu$s, while in 3D it is 0.2218
        $\mu$s.}
    \label{fig:single_cell_performance}
\end{figure}

To verify the computational complexity estimates, we isolate the computational cost of the core operator components
from parallel overhead and memory access patterns, we first perform a sequential benchmark on a single cell. This
experiment is designed to be compute-bound by repeating the evaluation many times over the same data, thus highlighting
the raw floating-point performance and algorithmic complexity of each method. To illustrate the performance of the key
components of the method, we benchmarked the computational cost of the Laplacian operator on a single cell in both 2D
and 3D. We compared three approaches: (i) a standard matrix-free evaluation exploiting tensor product factorization
(\texttt{FEEvaluation}), (ii) a version that skips the tensor product factorization (\texttt{FEPointEvaluation}), and
(iii) the ghost penalty evaluation on one face (\texttt{GhostPenalty}). The first two names refer to classes in
\texttt{deal.II}. The \texttt{FEEvaluation} and \texttt{FEPointEvaluation} approaches correspond to the matrix-free
evaluation techniques presented in Bergbauer et al.~\cite{bergbauer2024high}. Our work extends this framework by
introducing an efficient, tensor-product-based evaluation for the ghost penalty term, which was not addressed
in~\cite{bergbauer2024high}. The results, averaged~\cite{hoefler2015scientific} over 10,000 applications, are shown in
Figure~\ref{fig:single_cell_performance}. We show the relative timings, that is obtained by dividing the time for each
method by the time of FEEvaluation for $k=1$. In 2D that is 0.077 $\mu$s, while in 3D it is 0.2218 $\mu$s. To ensure a
fair comparison, the timings for \texttt{FEEvaluation} and \texttt{GhostPenalty}, which are vectorized over multiple
cells using SIMD instructions, were divided by the vectorization size. The experiment was performed in double
precision, and the vectorization size is 4 (using 256-bit AVX2 instructions). The \texttt{FEPointEvaluation} is
vectorized within one cell.

The computational cost of the ghost penalty scales similarly to \texttt{FEEvaluation}, as both methods involve tensor
product operations. \texttt{FEPointEvaluation} exhibits a higher computational cost, that is especially visible in 3D
focus. This behavior aligns with the expected computational complexity estimates, where tensor product factorization
reduces the overall cost of the operator application. Note the test is performed with minimal number of quadrature
points, while for a cut cell the number of quadrature points might be higher. Hence, the evalution on those cell might
be even more expensive.

We evaluate the computational efficiency of the matrix-free approach by measuring the throughput of a single
matrix-vector multiplication (vmult), measured in degrees of freedom (DoFs) per second. The number of DoFs corresponds
to the active degrees of freedom, i.e., those associated with cells that have a non-empty intersection with the domain
$\Omega$. Degrees of freedom on cells lying entirely outside $\Omega$ are excluded from the computation and timing.
Figure~\ref{fig:vmult_time_vs_degree} shows that the time per degree of freedom decreases with increasing polynomial
degree, indicating improved efficiency for higher-order elements.

\begin{figure}[!ht]
    \centering
    {\label{fig:vmult_time_complex_domain_2d}
        \begin{tikzpicture}
            \begin{axis}[
                    xlabel=Number of DoFs,
                    ylabel=DoFs/second,
                    xmode=log,
                    ymode=log,
                    legend pos=south east,
                    width=0.47\textwidth,
                    height=0.45\textwidth,
                    grid=both,
                    minor tick num=1,
                    log basis x=10,
                    log basis y=10,
                    title= 2D ] \addplot[color=blue, mark=o, mark size=3pt, line width=1pt] table [ x index=2, y
                        expr=\thisrowno{2}/\thisrowno{3}*1e3, header=has colnames ] {results/vmult_single_2D_d1.txt}; \addlegendentry{\(k=1\)}

                \addplot[color=red, mark=x, mark size=3pt, line width=1pt] table [
                        x index=2,
                        y expr=\thisrowno{2}/\thisrowno{3}*1e3,
                        header=has colnames
                    ] {results/vmult_single_2D_d2.txt};
                \addlegendentry{\(k=2\)}

                \addplot[color=green, mark=+, mark size=3pt, line width=1pt] table [
                        x index=2,
                        y expr=\thisrowno{2}/\thisrowno{3}*1e3,
                        header=has colnames]
                    {results/vmult_single_2D_d3.txt};
                \addlegendentry{\(k=3\)}
            \end{axis}
        \end{tikzpicture}
    }
    \hfill
    {\label{fig:vmult_time_vs_degree_3d}
        \begin{tikzpicture}
            \begin{axis}[
                    title= 3D,
                    xlabel=Number of DoFs,
                    ylabel=DoFs/second,
                    xmode=log,
                    ymode=log,
                    legend pos=south east,
                    width=0.47\textwidth,
                    height=0.45\textwidth,
                    grid=both,
                    log basis x=10,
                    log basis y=10,
                    ytick={1e6, 1e7, 1e8, 1e9},
                    ymin=4e7,
                    ymax=1.3e9,
                ]
                \addplot[color=blue, mark=o, mark size=3pt, line width=1pt] table [
                        x index=2,
                        y expr=\thisrowno{2}/\thisrowno{3}*1e3,
                        header=has colnames
                    ]
                    {results/vmult_single_3D_d1.txt};
                \addlegendentry{\(k=1\)}

                \addplot[color=red, mark=x, mark size=3pt, line width=1pt] table [
                        x index=2,
                        y expr=\thisrowno{2}/\thisrowno{3}*1e3,
                        header=has colnames
                    ] {results/vmult_single_3D_d2.txt};
                \addlegendentry{\(k=2\)}

                \addplot[color=green, mark=+, mark size=3pt, line width=1pt] table [
                        x index=2,
                        y expr=\thisrowno{2}/\thisrowno{3}*1e3,
                        header=has colnames]
                    {results/vmult_single_3D_d3.txt};
                \addlegendentry{\(k=3\)}
            \end{axis}
        \end{tikzpicture}
    }
    \caption{Problem with a single ball: Throughput of matrix-vector multiplication in degrees of freedom per second
        for different
        polynomial degrees. All throughput values are reported as DoFs/s, and only DoFs on active cells are
        included.}
    \label{fig:vmult_time_vs_degree}
\end{figure}

It is important to note that the method benefits from the tensor-product structure of the finite element space, which
can only be fully exploited on cells that are not intersected by the boundary. For cells that are intersected, the
integration requires special quadrature rules that break the tensor product structure, leading to higher computational
costs per element. However, as the number of intersected cells grows with $O(h^{d-1})$ while the total number of cells
grows with $O(h^{-d})$, their impact on the overall performance diminishes with mesh refinement.

The problem becomes more pronounced in case of more complex domains, where the ratio of intersected cells to the total
number of cells increases~\cite{bergbauer2024high}. To illustrate this, we consider a domain populated with multiple
balls, each with randomly chosen centers and a radius inversely proportional to the number of balls. Although the
domain may not necessarily be connected, as the balls might not intersect, this setup allows us to effectively
demonstrate the impact of the number of intersected cells on the overall performance.

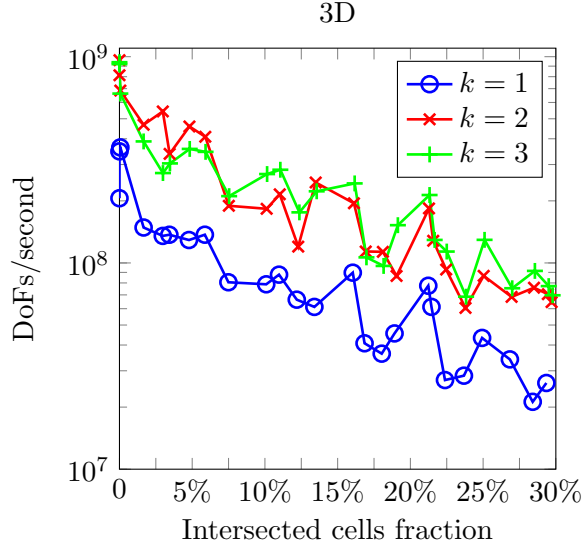
\begin{figure}[!ht]
    \centering
    {\label{fig:vmult_time_complex_domain_3d}
        \begin{tikzpicture}
            \begin{axis}[
                    title=3D,
                    xlabel=Intersected cells fraction,
                    ylabel=DoFs/second,
                    legend pos=north east,
                    ymode=log,
                    width=0.46\textwidth,
                    height=0.45\textwidth,
                    grid=both,
                    xmin=0,
                    xmax=0.30,
                    xtick={0, 0.025, 0.05, 0.075, 0.1, 0.125, 0.15, 0.175, 0.2, 0.225, 0.25,0.275, 0.30},
                    xticklabels={0, , 5\%, , 10\%, , 15\%, , 20\%, , 25\%, , 30\%},
                    ymin=1e7,
                    ymax=1.1e9,
                ]
                \addplot[color=blue, mark=o, mark size=3pt, line width=1pt] table [
                        x index=10,
                        y expr=\thisrowno{3}/\thisrowno{4},
                        header=has colnames,
                    ] {results/vmult_multiple_3D_d1_alt.txt};
                \addlegendentry{\(k=1\)}

                \addplot[color=red, mark=x, mark size=3pt, line width=1pt] table [
                        x index=10,
                        y expr=\thisrowno{3}/\thisrowno{4},
                        header=has colnames,
                    ] {results/vmult_multiple_3D_d2_alt.txt};
                \addlegendentry{\(k=2\)}

                \addplot[color=green, mark=+, mark size=3pt, line width=1pt] table [
                        x index=10,
                        y expr=\thisrowno{3}/\thisrowno{4},
                        header=has colnames]
                    {results/vmult_multiple_3D_d3_alt.txt};
                \addlegendentry{\(k=3\)}
            \end{axis}
        \end{tikzpicture}

    }
    \caption{
        Problem with multiple balls: Throughput of matrix-vector multiplication in DoFs per second for varying numbers
        of balls.}
    \label{fig:vmult_time_complex_domain}
\end{figure}

For a fixed 3D mesh, obtained by refining a single cell 5 times, we vary the number of balls and measure the vmult time
per unknown. Figure~\ref{fig:vmult_time_complex_domain} illustrates the throughput in DoFs per second for different
polynomial degrees as a function of the intersected cell fraction. Note that the number of degrees of freedom per cell
varies depending on the geometry of the domain, as the cells outside the domain are omitted in the computation. The
minimum and maximum number of degrees of freedom per cell for different polynomial degrees are shown in
Table~\ref{tab:dofs_per_cell}. The highest number of degrees of freedom per cell is achieved for the case of a single
ball in the domain, while the lowest is achieved for the case of a domain consisting of a maximum number of balls.

\begin{table}[!ht]
    \centering
    \caption{Minimum and maximum number of degrees of freedom per cell for different polynomial degrees in 3D,
        multiple ball problem.}
    \begin{tabular}{c |c c c c c c}
        Polynomial Degree &  & Min \#DoFs & Max \#DoFs \\ \hline 1 & & 251,428 & 2,146,689 \\ 2 & & 1,362,568 & 16,974,593 \\ 3
                          &  & 6,129,445  & 57,066,625 \\
    \end{tabular}

    \label{tab:dofs_per_cell}
\end{table}

We clearly observe that the throughput decreases with the intersected cell fraction, and this performance deterioration
is more pronounced for higher polynomial degrees. This is expected as the integration over the intersected cells
requires special quadrature rules that break the tensor product structure, leading to higher computational costs per
element, as illustrated in Figure~\ref{fig:single_cell_performance}.

\begin{figure}[!ht]
    \begin{minipage}[b][0.45\textwidth][c]{0.47\textwidth}
        \includegraphics[width=\linewidth]{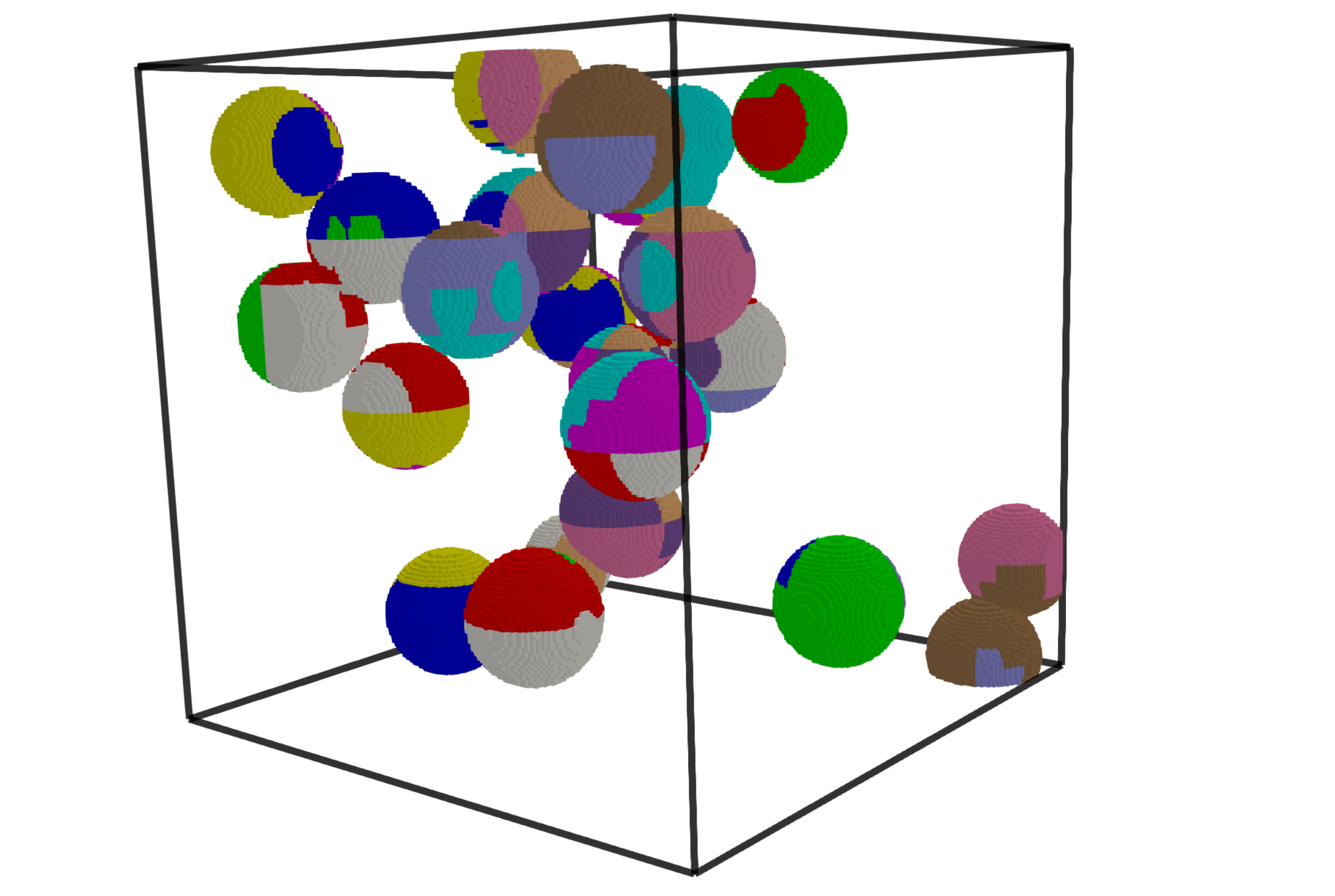}
    \end{minipage}\hfill
    \begin{tikzpicture}
        \begin{axis}[
                title= 3D,
                xlabel=Number of DoFs,
                ylabel=DoFs/second,
                xmode=log,
                ymode=log,
                legend pos=south east,
                width=0.47\textwidth,
                height=0.45\textwidth,
                grid=both,
                log basis x=10,
                log basis y=10,
                ytick={1e6, 1e7, 1e8, 1e9},
                ymin=4e7,
                ymax=1.3e9,
            ]
            \addplot[color=blue, mark=o, mark size=3pt, line width=1.2pt] table [
                    x index=3,
                    y expr=\thisrowno{3}/\thisrowno{4}*1e3,
                    header=has colnames
                ]
                {results/vmult_ref_mult_3D_d1.txt};
            \addlegendentry{\(k=1\)}

            \addplot[color=red, mark=x, mark size=3pt, line width=1.2pt] table [
                    x index=3,
                    y expr=\thisrowno{3}/\thisrowno{4}*1e3,
                    header=has colnames
                ] {results/vmult_ref_mult_3D_d2.txt};
            \addlegendentry{\(k=2\)}

            \addplot[color=green, mark=+, mark size=3pt, line width=1.2pt] table [
                    x index=3,
                    y expr=\thisrowno{3}/\thisrowno{4}*1e3,
                    header=has colnames
                ] {results/vmult_ref_mult_3D_d3.txt};
            \addlegendentry{\(k=3\)}

            \addplot[color=blue, mark=o, mark size=2pt, line width=1pt, dashed] table [
                    x index=2,
                    y expr=\thisrowno{2}/\thisrowno{3}*1e3,
                    header=has colnames
                ]
                {results/vmult_single_3D_d1.txt};

            \addplot[color=red, mark=x, mark size=2pt, line width=1pt, dashed] table [
                    x index=2,
                    y expr=\thisrowno{2}/\thisrowno{3}*1e3,
                    header=has colnames
                ] {results/vmult_single_3D_d2.txt};

            \addplot[color=green, mark=+, mark size=2pt, line width=1pt, dashed] table [
                    x index=2,
                    y expr=\thisrowno{2}/\thisrowno{3}*1e3,
                    header=has colnames
                ] {results/vmult_single_3D_d3.txt};

        \end{axis}
    \end{tikzpicture}
    \caption{Left: Test problem domain with 25 balls, 29\% intersected cells. Colors indicate the MPI rank responsible
        for the
        cell. Black lines form the outline of the mesh. Right: Throughput of matrix-vector multiplication in DoF per
        second for a problem with 25 balls on
        subsequently refined meshes. The dashed lines show the performance for a single ball for comparison.
        Throughput is measured as DoFs/s, counting only active
        degrees of freedom.
    }
    \label{fig:vmult_time_vs_refinement}
\end{figure}

To show the diminishing impact of the intersected cells on the overall performance, we repeat the test for a fixed
number of balls (25), resulting in the fraction of intersected cells decreasing with mesh refinement from 60\% to 15\%.
The problem domain for 25 balls is illustrated in the left part of Figure~\ref{fig:vmult_time_vs_refinement}. We
measure throughput for selected polynomial degrees on subsequently refined meshes. The right panel of
Figure~\ref{fig:vmult_time_vs_refinement} illustrates the performance of the method with increasing problem size.

More insights can be gained by breaking down the execution time into the time spent on each part of the algorithm. We
examine a 3D problem from the previous experiments, using the finest mesh, $k=3$, and 26\% cut cells, resulting in a
total of 6,129,445 DoFs.

\begin{figure}[!ht]
    \centering
    \begin{tikzpicture}
        \begin{axis}[
                xbar stacked,
                xmin=0, xmax=100,
                xlabel={Percentage of Total vmult Time},
                width=0.8\textwidth,
                height=0.3\textwidth,
                ytick=\empty,
                axis lines=none
            ]
            \addplot[fill=blue!40, draw=black, line width=0.5pt] coordinates {(6.50,0)};
            \node at (axis cs:8.50/2,0.25) [blue!80!black]
            {6.50\%};
            \node at (axis cs:10.50/2,-0.25) [blue!80!black] {Interior};

            \addplot[fill=red!40, draw=black, line width=0.5pt] coordinates {(65.5,0)};
            \node at (axis cs:6.50 + 65.5/2,0.25) [red!80!black]
            {65.5\%};
            \node at (axis cs:6.50 + 65.5/2,-0.25) [red!80!black]
            {Intersected};

            \addplot[fill=green!40, draw=black, line width=0.5pt] coordinates {(14.6,0)};
            \node at (axis cs:6.50 + 65.5 + 14.6/2,0.25)
            [green!80!black] {14.6\%};
            \node at (axis cs:6.50 + 65.5 + 12.6/2,-0.25)
            [green!80!black] {Ghost Penalty};

            \addplot[fill=gray!40, draw=black, line width=0.5pt] coordinates {(12.3,0)};
            \node at (axis cs:6.50 + 65.5 + 14.6 + 12.3/2,0.25)
            [gray!80!black] {12.3\%};
            \node at (axis cs:6.50 + 65.5 + 14.6 + 12.3/2,-0.25)
            [gray!80!black] {MPI};

        \end{axis}
    \end{tikzpicture}
    \caption{Breakdown of vmult time into different components for 3D problem with  $k=3$, 17\% cut  cells, total
        number of DoFs: 6,129,445 }
    \label{fig:vmult_time_breakdown}
\end{figure}
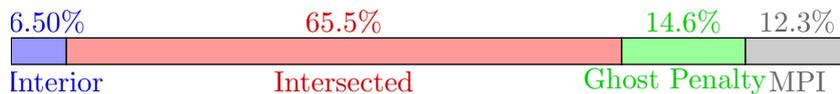

Figure~\ref{fig:vmult_time_breakdown} shows the percentage of time spent on computations involving cut cells, interior
cells, ghost penalty terms, and the time spent on MPI communication. The timings are aggregated over all 32 MPI ranks.
The 'MPI' category represents time spent in communication and synchronization, and therefore implicitly includes the
effects of load imbalance. The dominance of the 'Intersected' cells category is consistent with the single-cell
benchmarks in Figure~\ref{fig:single_cell_performance}, which show that evaluation on cut cells (which cannot use
tensor-product factorization for quadrature) is significantly more expensive than on interior cells. The cut cells
consume the most time (65.5\%), followed by the ghost penalty (14.6\%), MPI communication (12.3\%) and interior cells
(6.5\%). The interior cell computations are more efficient due to the exploitation of the tensor-product structure that
is not possible on the cut cells. Computing the ghost penalty takes a fraction of the time compared to the cut cells.
The MPI communication time is relatively low, given how simple the load balancing algorithm is.

\section{Conclusion}
\label{sec:conclusion}

Conclusion was rewritten.

We have presented an efficient matrix-free approach for implementing and evaluating ghost penalty stabilization within
the Cut Finite Element Method framework. Our contribution lies in recognizing and exploiting the tensor-product
structure inherent in the ghost penalty operator, which allows us to decompose the evaluation into a sequence of
one-dimensional operations. This factorization reduces the computational complexity to $O(k^{d+1})$ for elements of
degree $k$ in $d$ dimensions, making high-order stabilization computationally feasible. While previous works have
obtained similar complexity benefits~\cite{bergbauer2024high}, our approach is novel in its application to the ghost
penalty term.

The presented method completely avoids the assembly and storage of global matrices, instead relying on precomputed
one-dimensional operators and geometrical data. This approach provides substantial memory savings, particularly for
high-order methods where storage requirements for assembled matrices grow rapidly with polynomial degree. Our numerical
experiments demonstrate that the method achieves optimal convergence rates while maintaining high computational
efficiency, even as the problem size increases.

A notable advantage of our approach is that the computational cost of the ghost penalty evaluation scales similarly to
the matrix-free evaluation of standard finite element operators using tensor product factorization. As shown in our
performance analysis, while cut cells remain the most expensive component of the computation, the ghost penalty
evaluation contributes only a modest fraction to the overall execution time. Furthermore, the implementation within the
\texttt{deal.II} library allows further development of the method.

\paragraph*{Code availability:} \url{https://github.com/mwichro/TensorGhostPenalty}
\section*{Acknowledgments}
During the preparation of this work, the author used Google Gemini, ChatGPT, and Claude Sonnet in order to improve the
clarity and readability of the manuscript. After using these tools, the author reviewed and edited the content as
needed and takes full responsibility for the content of the publication.

\bibliographystyle{siam}
\bibliography{literature}

\end{document}